\documentclass[11pt]{amsart}
\usepackage{txfonts}
\usepackage{mathrsfs}
\usepackage{amsmath}
\usepackage[pagewise]{lineno}
\allowdisplaybreaks[4]
\usepackage{}
\usepackage{float}
\usepackage{xypic}
\usepackage{amsfonts}
\usepackage{amssymb}
\usepackage{bbm}
\usepackage{color,xcolor}

\usepackage[colorlinks,
            linkcolor=blue,
           anchorcolor=blue,
            citecolor=blue
            ]{hyperref}



\newtheorem{thm}{Theorem}[section]
\newtheorem{cor}[thm]{Corollary}
\newtheorem{lem}[thm]{Lemma}
\newtheorem{exa}[thm]{Example}
\newtheorem{prop}[thm]{Proposition}
\newtheorem{rem}[thm]{Remark}
\numberwithin{equation}{section}

\newcommand{\be}{\begin{equation}}
\newcommand{\ee}{\end{equation}}
\newcommand{\bes}{\begin{eqnarray}}
\newcommand{\ees}{\end{eqnarray}}
\newcommand{\bess}{\begin{eqnarray*}}
\newcommand{\eess}{\end{eqnarray*}}

\newcommand{\bali}{\begin{align}}
\newcommand{\eali}{\end{align}}

\def\kk{\mathbbm{k}}

\begin{document}
\title[Invariants from the Sweedler power maps on integrals]{Invariants from the Sweedler power maps on integrals}
\author{Zhihua Wang}
\address{Z. Wang\newline Department of Mathematics, Taizhou University,
Taizhou 225300, China}
\email{mailzhihua@126.com}
\author{Gongxiang Liu}
\address{G. Liu\newline Department of Mathematics, Nanjing University,
Nanjing 210093, China}
\email{gxliu@nju.edu.cn}
\author{Libin Li}
\address{L. Li\newline School of Mathematical Science, Yangzhou University, Yangzhou 225002, China}
\email{lbli@yzu.edu.cn}
\date{}

\begin{abstract}For a finite-dimensional Hopf algebra $A$ with a nonzero left integral $\Lambda$, we investigate a relationship between $P_n(\Lambda)$ and $P_n^J(\Lambda)$, where $P_n$ and $P_n^J$ are respectively the $n$-th Sweedler power maps of $A$ and the twisted Hopf algebra $A^J$. We use this relation to give several invariants of the representation category Rep$(A)$ considered as a tensor category. As applications, we distinguish the representation categories of 12-dimensional
pointed nonsemisimple Hopf algebras. Also, these invariants are sufficient to distinguish the representation categories Rep$(K_8)$, Rep$(\kk Q_8)$ and Rep$(\kk D_4)$,
although they have been completely distinguished by their Frobenius-Schur indicators.
We further reveal a relationship between the right integrals $\lambda$ in $A^*$ and $\lambda^J$ in $(A^J)^*$. This can be used to give a uniform proof of the remarkable result which says that the $n$-th indicator $\nu_n(A)$ is a gauge invariant of $A$ for any $n\in \mathbb{Z}$. We also use the expression for $\lambda^J$ to give an alternative proof of the known result that the Killing form of the Hopf algebra $A$ is invariant under twisting. As a result, the dimension of the Killing radical of $A$ is a gauge invariant of $A$.
\end{abstract}
\maketitle
\section{\bf Introduction}
It is very interesting to know when two finite-dimensional Hopf algebras are gauge equivalent, in the sense that the categories of their finite-dimensional representations are equivalent as tensor categories. It can be seen from  \cite[Theorem 2.2]{NS1} that two finite-dimensional Hopf algebras $A$ and $A'$ are gauge equivalent if and only if there exists a twist $J$ on $A$ such that $A'$ and the twisted Hopf algebra $A^J$ are isomorphic as bialgebras. However, to determine the existence or nonexistence of such twists and isomorphisms is still highly non-trivial.

A quantity $\kappa(A)$ defined for any finite-dimensional
Hopf algebra $A$ is called a gauge invariant of $A$ if $\kappa(A')=\kappa(A)$
for any Hopf algebra $A'$ gauge equivalent to $A$. For instance, the dimension, the quasi-exponent (see \cite{EG}) and the $n$-th indicator (see \cite{Sh}) of any finite-dimensional Hopf algebra are gauge invariants. We refer to \cite{KMM,KMN,KSZ,CS,WW} for more details about gauge invariants of Hopf algebras. These invariants are used as powerful tools to distinguish representation categories of Hopf algebras. However, more general but computable invariants are still needed to be discovered.


The present paper is specifically concerned with those invariants raised from the images of the Sweedler power maps on integrals.
For any finite-dimensional Hopf algebra $A$ with a nonzero left integral $\Lambda$ and a normalized twist $J$, the present paper is dedicated to a study of the relationship between $P_n(\Lambda)$ and $P_n^J(\Lambda)$, where $P_n$ and $P_n^J$ are respectively the $n$-th Sweedler power maps of the Hopf algebra $A$ and the twisted Hopf algebra $A^J$. We show that
$$P^J_n(\Lambda)=TP_n(\Lambda)\ \text{for}\ n\in\mathbb{Z},$$ where $T$ is an invertible element of $A$ associated to the twist $J$ (see Theorem \ref{cor03}).
There are some invariants that can be derived from this equality. For instance, we show that the homogeneous polynomials vanishing on some $P_n(\Lambda)$ are invariants of the representation category Rep$(A)$ (see Theorem \ref{p1}). As an application, we show that the representation categories of 12-dimensional pointed nonsemisimple Hopf algebras classified in \cite{Natale} are mutually inequivalent as tensor categories.

If, moreover, $A$ is unimodular, then $T=1$ and hence $P^J_n(\Lambda)=P_n(\Lambda)$ for $n\in\mathbb{Z}$, namely, $P_n(\Lambda)$ is invariant under twisting. Even in this case, the (not necessarily homogeneous) polynomials vanishing on some $P_n(\Lambda)$ are in general not invariants of Rep$(A)$. However, this kind of polynomials can also be used to distinguish the representation categories of Hopf algebras in a certain sense (see Proposition \ref{p2}). If $A$ is semisimple and $\Lambda$ is idempotent, then the polynomials vanishing on some $P_n(\Lambda)$ are indeed invariants of the representation category Rep$(A)$. We use this approach to distinguish the representation categories Rep$(K_8)$, Rep$(\kk Q_8)$ and Rep$(\kk D_4)$ whereas they have the same fusion rules, we refer to \cite{TY} and \cite{NS1} for other approaches to distinguish them.

Note that these $P_n(\Lambda)$ are dependent on the choice of $\Lambda$. The values that characters of finite-dimensional $A$-modules taking on them are not gauge invariants. In view of this, we turn to consider the ratio $$\chi_{V}(P_n(\Lambda)):\chi_{W}(P_m(\Lambda))$$ for finite-dimensional $A$-modules $V$ and $W$.
When  $A$ is unimodular and $\chi_{W}(P_m(\Lambda))\neq0$, we show that this ratio is a gauge invariant of $A$. This invariant  can be regarded as a slightly generalization of the $n$-th FS-indicator defined for any semisimple Hopf algebra. Indeed, if $A$ is semisimple and $W$ is chosen to be the trivial $A$-module $\kk$, then the above ratio becomes $\chi_{V}(P_n(\Lambda/\varepsilon(\Lambda)))$, which is the $n$-th FS-indicator of $V$ (see \cite{KSZ,LM}).



We further describe a relationship between the right integrals $\lambda$ in $A^*$ and $\lambda^J$ in $(A^J)^*$. We show that
$$\lambda^J=\lambda\leftharpoonup S^2(R^{-1})S(Q_J^{-1})Q_J,$$ where $R$ and $Q$ are invertible elements of $A$ associated to the twist $J$. This generalizes the result of \cite[Theorem 3.4]{AEGN} for a unimodular Hopf algebra.
We use this formula to give a uniform proof of the remarkable result that the $n$-th indicator $\nu_n(A)$ of $A$ is a gauge invariant for any $n\in \mathbb{Z}$. We also use this formula to give an alternative proof of the known result that the Killing form of the Hopf algebra $A$ is invariant under twisting. This gives rise to a result that the dimension of the Killing radical of $A$ is a gauge invariant of $A$.

The paper is organized as follows. In Section 2, we recall some basic facts about Hopf algebras and twisted Hopf algebras. In Section 3, we mainly investigate a relationship  between $P^J_n(\Lambda)$ and $P_n(\Lambda)$ for any $n\in\mathbb{Z}$.  The importance of the relationship is illustrated in Section 4: Some invariants of representation categories of finite-dimensional Hopf algebras can be deduced from the relationship and then they are used to distinguish the representation categories of some Hopf algebras. In Section 5, we show that a right integral $\lambda^J$ in $(A^J)^*$ can be expressed in terms of a right integral $\lambda$ in $A^*$. This enables us to give a uniform proof of the result that the $n$-th indicator $\nu_n(A)$ for $n\in\mathbb{Z}$ is a gauge invariant of $A$. This also enables us to give an alternative proof of the known result that the Killing form of the Hopf algebra $A$ is invariant under twisting. 

\section{\bf Preliminaries}
Throughout the paper we work over a fixed base field $\kk$. Hopf algebras, vector spaces,
and tensor products are understood to be over $\kk$. For a finite-dimensional Hopf algebra $A$ over the field $\kk$, we use Sweedler notation in the form
$\Delta(a)=a_{(1)}\otimes a_{(2)}$ for comultiplication $\Delta$ of $A$ and write $S$ for the antipode of $A$. The inverse of $S$ under composition is denoted by $S^{-1}$.
We fix a left integral $\Lambda$ in $A$ and a right integral $\lambda$  in $A^*$ such that $\lambda(\Lambda)=1$. We denote $\alpha$ to be the distinguished group-like element in $A^*$ defined by $\Lambda a=\alpha(a)\Lambda$ for $a\in A$. The category of finite-dimensional left $A$-modules is denoted by Rep$(A)$, which is a rigid tensor category.

The following equalities are of fundamental importance in this paper. 
The first equality can be found in \cite[Lemma 1.2]{LR}:
\begin{equation}\label{equ004}S(a)\Lambda_{(1)}\otimes \Lambda_{(2)}=\Lambda_{(1)}\otimes a\Lambda_{(2)}\ \text{for}\ a\in A.\end{equation}
The second can be found in \cite{KMN}:
\begin{equation}\label{equ005}\Lambda_{(1)}a\otimes \Lambda_{(2)}=\Lambda_{(1)}\otimes \Lambda_{(2)}\alpha(a_{(1)})S(a_{(2)})\ \text{for}\ a\in A.\end{equation}
For $a,b\in A$, the integral $\lambda$ satisfies the following properties:
\begin{gather}
\label{equ4004}\lambda(ab)=\lambda(\alpha(b_{(1)})S^2(b_{(2)})a),\\
\label{equ4004444}\lambda(ab_{(1)})b_{(2)}=\lambda(a_{(1)}b)S^{-1}(a_{(2)}),
\end{gather}
where the former follows from \cite[Theorem 3(a)]{Rad} and the latter follows from \cite[Remark 3.2]{AEGN}.
For any $\kk$-linear map $f:A\rightarrow A$, the trace of $f$ can be described by Radford's trace formula (see \cite[Theorem 2]{Rad}): \begin{equation}\label{equ400444}\text{tr}(f)=\lambda(S(\Lambda_{(2)})f(\Lambda_{(1)})).\end{equation}

Recall from \cite{LMS} that the $n$-th Sweedler power map $P_n:A\rightarrow A$ is the $n$-th convolution power of the identity map $id_A$ of $A$. Namely,
\begin{equation*}\label{equ0}P_n(a)=\left\{
           \begin{array}{ll}
             a_{(1)}\cdots a_{(n)}, &  n\geq1; \\
             \varepsilon(a), & n=0; \\
             S(a_{(-n)}\cdots a_{(1)}), &  n\leq-1.
           \end{array}
         \right.
\end{equation*}
Note that $P_1(a)=a$ and $P_{-1}(a)=S(a)$ for $a\in A$.


Recall that a normalized twist for a finite-dimensional Hopf algebra $A$ is an invertible element $J\in A\otimes A$ which satisfies $(\varepsilon\otimes id)(J)=(id\otimes\varepsilon)(J)=1$ and
\begin{equation}\label{qq1}(\Delta\otimes id)(J)(J\otimes1)=(id\otimes\Delta)(J)(1\otimes J).\end{equation}
We write $J=J^{(1)}\otimes J^{(2)}$ and $J^{-1}=J^{-(1)}\otimes J^{-(2)}$, where the summation is understood. We also write $J_{21}=J^{(2)}\otimes J^{(1)}$.

Given a normalized twist $J$ for $A$ one can define a new Hopf algebra $A^{J}$ with the same algebra structure and counit as $A$, for which the comultiplication $\Delta^J$ and antipode $S^J$ are given respectively by
$$\Delta^J(a)=J^{-1}\Delta(a)J,\ \ \ \
S^J(a)=Q^{-1}_JS(a)Q_J\ \text{for}\ a\in A,$$
where $Q_J=S(J^{(1)})J^{(2)}$, which is an invertible element of $A$ with the inverse $Q_J^{-1}=J^{-(1)}S(J^{-(2)})$.

The element $Q_J$ satisfies the following identity (see \cite[Eq (5)]{AEGN}):
\begin{equation}\label{equ002}
\Delta(Q_J)=(S\otimes S)(J^{-1}_{21})(Q_J\otimes Q_J)J^{-1}.
\end{equation}
This implies that
\begin{equation}\label{equu00003}
\Delta(Q_J^{-1})=J(Q_J^{-1}\otimes Q_J^{-1})(S\otimes S)(J_{21}).
\end{equation}
The following three equalities follow from \cite[Lemma 2.4]{AEGN}:
\begin{gather}
\label{equu0003}(1\otimes Q_J^{-1})(S\otimes S)(J_{21})=(1\otimes J^{-(1)})\Delta(S(J^{-(2)})),\\
\label{equ003}J_{(1)}^{(1)}\otimes S(J^{(1)}_{(2)})J^{(2)}=J^{-(1)}\otimes S(J^{-(2)})Q_J,\\
\label{equ0033}J^{-(1)}_{(1)}\otimes J^{-(1)}_{(2)}S(J^{-(2)})=J^{(1)}\otimes J^{(2)}Q_J^{-1}.
\end{gather}

We define $\Delta^{(1)}=id_{A}$ and $\Delta^{(n+1)}=(id_{A}\otimes \Delta^{(n)})\circ\Delta$ for all $n\geq1$. Note that a nonzero left integral $\Lambda$ of $A$ is also a left integral of $A^{J}$. We denote $$\Delta^{(n)}(\Lambda)=\Lambda_{(1)}\otimes\cdots\otimes\Lambda_{(n)}\ \text{and}\ (\Delta^{J})^{(n)}(\Lambda)=\Lambda_{\langle1\rangle}\otimes\cdots\otimes\Lambda_{\langle n\rangle}$$ to distinguish between $\Delta^{(n)}(\Lambda)$ and $(\Delta^{J})^{(n)}(\Lambda).$ The $n$-th Sweedler power map of $A^J$  is denoted by $P_n^J$. Then $P_n^J(a)$ can be written as
\begin{equation*}\label{equ0}P^J_n(a)=\left\{
           \begin{array}{ll}
             a_{\langle1\rangle}\cdots a_{\langle n\rangle}, &  n\geq1; \\
             \varepsilon(a), & n=0; \\
             S^J(a_{\langle-n\rangle}\cdots a_{\langle1\rangle}), &  n\leq-1.
           \end{array}
         \right.
\end{equation*}



\section{\bf The Sweedler power maps on integrals}
In this section, we first give a property of $P_n(\Lambda)$  for a left integral $\Lambda$ in $A$. We then use this property to describe a relationship between $P_n^J(\Lambda)$ and $P_n(\Lambda)$ for a normalized twist $J$ of $A$.

\begin{prop}\label{lem1}
Let $A$ be a finite-dimensional Hopf algebra over the field $\kk$ with a nonzero left integral $\Lambda$.  For any $n\in\mathbb{Z}$ and $a\in A$, we have $aP_n(\Lambda)=P_n(\Lambda)a^{\dag}$, where $a^{\dag}=a_{(1)}(\alpha\circ S^{-1})(a_{(2)}).$
Moreover, if $A$ is unimodular, then $P_n(\Lambda)$ is a central element of $A$.
\end{prop}
\proof 
Obviously, the desired result holds for the cases $n=0,\pm1$. For the case $n>1$, we have
\begin{align*}
aP_n(\Lambda)&=a\Lambda_{(1)}P_{n-1}(\Lambda_{(2)})\\
&=\Lambda_{(1)}P_{n-1}(S^{-1}(a)\Lambda_{(2)})\ \  \text{by}\ (\ref{equ004})\\
&=\Lambda_{(1)}(S^{-1}(a))_{(1)}\Lambda_{(2)}(S^{-1}(a))_{(2)}\cdots\Lambda_{(n-1)}(S^{-1}(a))_{(n-1)}\Lambda_{(n)}\\
&=P_{n-1}(\Lambda_{(1)}S^{-1}(a))\Lambda_{(2)}\\
&=P_{n-1}(\Lambda_{(1)})\Lambda_{(2)}a_{(1)}\alpha(S^{-1}(a_{(2)}))\ \  \text{by}\ (\ref{equ005})\\
&=P_n(\Lambda)a^{\dag}.
\end{align*}
For the case $n<-1$, we have
\begin{align*}
aP_n(\Lambda)&=S(\Lambda_{(1)}S^{-1}(a))P_{n+1}(\Lambda_{(2)})\\
&=S(\Lambda_{(1)})P_{n+1}(\Lambda_{(2)}a^{\dag})\ \   \text{by}\ (\ref{equ005})\\
&=S(\Lambda_{(1)})S(\Lambda_{(-n)}a^{\dag}_{(-n-1)}\cdots\Lambda_{(2)}a^{\dag}_{(1)})\\
&=S(a^{\dag}_{(-n-1)}\Lambda_{(-n-1)}\cdots a^{\dag}_{(1)}\Lambda_{(1)})S(\Lambda_{(-n)})\\
&=P_{n+1}(a^{\dag}\Lambda_{(1)})S(\Lambda_{(2)})\\
&=P_{n+1}(\Lambda_{(1)})S(S^{-1}(a^{\dag})\Lambda_{(2)})\ \   \text{by}\ (\ref{equ004})\\
&=S(\Lambda_{(-n-1)}\cdots\Lambda_{(1)})S(\Lambda_{(-n)})a^{\dag}\\
&=P_n(\Lambda)a^{\dag}.
\end{align*}
If $A$ is unimodular, then $\alpha=\varepsilon$ and hence $a^{\dag}=a$. In this case, $P_n(\Lambda)$ is a central element of $A$.
\qed

\begin{rem}For any $a\in A$, we define a family of elements $a_k\in A$ recursively by
\begin{equation}\label{er}a_1=a_{(1)}(\alpha\circ S^{-1})(a_{(2)})\ \text{and}\ a_{k+1}=(a_{k})_{(1)}(\alpha\circ S^{-1})((a_{k})_{(2)})\ \text{for}\ k\geq1.\end{equation} By induction on $k$, we may see that $a(P_n(\Lambda))^k=(P_n(\Lambda))^ka_k$  for any $n\in\mathbb{Z}.$
This further implies that $$(aP_n(\Lambda))^k=(P_n(\Lambda))^ka_k\cdots a_2a_1.$$
Moreover, for positive integers $k_1,\cdots,k_s$ and integers $n_1,\cdots,n_s,$ we have
\begin{equation}\label{eer}(aP_{n_1}(\Lambda))^{k_1}\cdots (aP_{n_s}(\Lambda))^{k_s}
=(P_{n_1}(\Lambda))^{k_1}\cdots(P_{n_s}(\Lambda))^{k_s}a_{k_1+\cdots+k_s}\cdots a_2a_1.
\end{equation}
\end{rem}

To describe a relationship between $P_n^J(\Lambda)$ and $P_n(\Lambda)$, we need the following lemma:

\begin{lem}\label{prop1} Let $A$ be a finite-dimensional Hopf algebra over the field $\kk$ with a normalized twist $J$ and a nonzero left integral $\Lambda$. Denote by $R=\alpha(J^{-(1)})J^{-(2)}$ and $T=J^{-(1)}\alpha(J^{-(2)}).$ For any integer number $n\geq2$, we have
\begin{enumerate}
  \item $\Lambda_{\langle1\rangle}\otimes \Lambda_{\langle2\rangle}\Lambda_{\langle3\rangle}\cdots\Lambda_{\langle n\rangle}$
  $$=Q_J^{-1}\Lambda_{(1)}J^{(1)}\otimes J^{-(1)}\Lambda_{(2)}P_{n-2}\bigg(J^{(2)}J^{-(2)}\Lambda_{(3)}\bigg)S(R)Q_J.$$
  \item $\Lambda_{\langle n\rangle}\otimes \Lambda_{\langle n-1\rangle}\Lambda_{\langle n-2\rangle}\cdots \Lambda_{\langle 1\rangle}$
      $$=S^{-1}(Q_J^{-1})\Lambda_{(3)}J^{(2)}\otimes J^{-(2)}\Lambda_{(2)}\bigg(S^{-1}\circ P_{2-n}\bigg)\bigg(J^{(1)}J^{-(1)}\Lambda_{(1)}\bigg)S^{-1}(T)S^{-1}(Q_J).$$
\end{enumerate}
\end{lem}
\proof We only give a proof of Part (1) and the proof of Part (2) is similar to that of Part (1).
We proceed by induction on $n$. For the case $n=2$,
\begin{align*}
\Lambda_{\langle1\rangle}\otimes\Lambda_{\langle2\rangle}&=\Delta^J(\Lambda)=J^{-1}\Delta(\Lambda)J=J^{-(1)}\Lambda_{(1)}J^{(1)}\otimes J^{-(2)}\Lambda_{(2)}J^{(2)}\\
&=J^{-(1)}S(J^{-(2)})\Lambda_{(1)}J^{(1)}\otimes \Lambda_{(2)}J^{(2)}\ \ \ \ \text{by}\ (\ref{equ004})\\
&=J^{-(1)}S(J^{-(2)})\Lambda_{(1)}\otimes \Lambda_{(2)}\alpha(J^{(1)}_{(1)})S(J^{(1)}_{(2)})J^{(2)} \ \ \ \text{by}\ (\ref{equ005})\\
&=Q_J^{-1}\Lambda_{(1)}\otimes \Lambda_{(2)}S(R)Q_J\ \ \ \ \text{by}\ (\ref{equ003})
\end{align*}
So the equality holds for the case $n=2$.
Suppose the identity in question holds for the case $n$. Namely,
\begin{equation}\label{equ001}\Lambda_{\langle1\rangle}\otimes\Lambda_{\langle2\rangle}\Lambda_{\langle3\rangle}\cdots\Lambda_{\langle n\rangle}=Q_J^{-1}\Lambda_{(1)}J^{(1)}\otimes J^{-(1)}\Lambda_{(2)}P_{n-2}\bigg(J^{(2)}J^{-(2)}\Lambda_{(3)}\bigg)S(R)Q_J.
\end{equation}
To prove the case $n+1$, we set $\widetilde{J}^{(1)}\otimes \widetilde{J}^{(2)}=J^{(1)}\otimes J^{(2)}$, $\widetilde{J}^{-(1)}\otimes \widetilde{J}^{-(2)}=J^{-(1)}\otimes J^{-(2)}$ and apply $\Delta^J\otimes id$ to both sides of (\ref{equ001}). It follows that
\begin{align*}
&\ \ \ \ \Lambda_{\langle1\rangle}\otimes\Lambda_{\langle2\rangle}\otimes\Lambda_{\langle3\rangle}\cdots\Lambda_{\langle n+1\rangle}\\
&=J^{-1}\Delta(Q_J^{-1})\Delta(\Lambda_{(1)})\Delta(J^{(1)})J\otimes J^{-(1)}\Lambda_{(2)}P_{n-2}(J^{(2)}J^{-(2)}\Lambda_{(3)})S(R)Q_J\\
&=(Q_J^{-1}\otimes1)\bigg((1\otimes Q_J^{-1})(S\otimes S)(J_{21})\bigg)\Delta(\Lambda_{(1)})\Delta(J^{(1)})J\otimes J^{-(1)}\Lambda_{(2)}P_{n-2}(J^{(2)}J^{-(2)}\Lambda_{(3)})S(R)Q_J\ \ \ \text{by}\ (\ref{equu00003})\\
&=(Q_J^{-1}\otimes1)(1\otimes \widetilde{J}^{-(1)})\Delta(S({\widetilde{J}}^{-(2)}))\Delta(\Lambda_{(1)})\Delta(J^{(1)})J\otimes J^{-(1)}\Lambda_{(2)}P_{n-2}(J^{(2)}J^{-(2)}\Lambda_{(3)})S(R)Q_J\ \ \ \text{by}\ (\ref{equu0003}) \\
&=(Q_J^{-1}\otimes \widetilde{J}^{-(1)})\Delta(S({\widetilde{J}}^{-(2)})\Lambda_{(1)})\Delta(J^{(1)})J\otimes J^{-(1)}\Lambda_{(2)}P_{n-2}(J^{(2)}J^{-(2)}\Lambda_{(3)})S(R)Q_J\\
&=(Q_J^{-1}\otimes \widetilde{J}^{-(1)})\Delta(\Lambda_{(1)})\Delta(J^{(1)})J\otimes J^{-(1)}{\widetilde{J}}^{-(2)}_{(1)}\Lambda_{(2)}P_{n-2}(J^{(2)}J^{-(2)}{\widetilde{J}}^{-(2)}_{(2)}\Lambda_{(3)})S(R)Q_J\ \ \ \text{by}\ (\ref{equ004})\\
&=(Q_J^{-1}\otimes \widetilde{J}^{-(1)})\Delta(\Lambda_{(1)})(\widetilde{J}^{(1)}\otimes \widetilde{J}^{(2)}_{(1)}J^{(1)})\otimes J^{-(1)}{\widetilde{J}}^{-(2)}_{(1)}\Lambda_{(2)}P_{n-2}(\widetilde{J}^{(2)}_{(2)}J^{(2)}J^{-(2)}{\widetilde{J}}^{-(2)}_{(2)}\Lambda_{(3)})S(R)Q_J\ \ \ \text{by}\ (\ref{qq1})\\
&=Q_J^{-1}\Lambda_{(1)}\widetilde{J}^{(1)}\otimes \widetilde{J}^{-(1)}\Lambda_{(2)}\widetilde{J}^{(2)}_{(1)}J^{(1)}\otimes J^{-(1)}\widetilde{J}^{-(2)}_{(1)}\Lambda_{(3)}P_{n-2}(\widetilde{J}^{(2)}_{(2)}J^{(2)}J^{-(2)}{\widetilde{J}}^{-(2)}_{(2)}\Lambda_{(4)})S(R)Q_J
\end{align*}
It follows that
\begin{align*}
&\ \ \ \ \Lambda_{\langle1\rangle}\otimes\Lambda_{\langle2\rangle}\Lambda_{\langle3\rangle}\cdots\Lambda_{\langle n+1\rangle}\\
&=Q_J^{-1}\Lambda_{(1)}\widetilde{J}^{(1)}\otimes \widetilde{J}^{-(1)}\Lambda_{(2)}\widetilde{J}^{(2)}_{(1)}J^{(1)} J^{-(1)}\widetilde{J}^{-(2)}_{(1)}\Lambda_{(3)}P_{n-2}(\widetilde{J}^{(2)}_{(2)}J^{(2)}J^{-(2)}{\widetilde{J}}^{-(2)}_{(2)}\Lambda_{(4)})S(R)Q_J\\
&=Q_J^{-1}\Lambda_{(1)}\widetilde{J}^{(1)}\otimes \widetilde{J}^{-(1)}\Lambda_{(2)}\widetilde{J}^{(2)}_{(1)}\widetilde{J}^{-(2)}_{(1)}\Lambda_{(3)}P_{n-2}(\widetilde{J}^{(2)}_{(2)}{\widetilde{J}}^{-(2)}_{(2)}\Lambda_{(4)})S(R)Q_J\\
&=Q_J^{-1}\Lambda_{(1)}\widetilde{J}^{(1)}\otimes \widetilde{J}^{-(1)}\Lambda_{(2)}P_{n-1}(\widetilde{J}^{(2)}{\widetilde{J}}^{-(2)}\Lambda_{(3)})S(R)Q_J\\
&=Q_J^{-1}\Lambda_{(1)}J^{(1)}\otimes J^{-(1)}\Lambda_{(2)}P_{n-1}(J^{(2)}J^{-(2)}\Lambda_{(3)})S(R)Q_J.
\end{align*}The equality of Part (1) is now proved by induction on $n$.
\qed

We have the following relationship between $P_n^J(\Lambda)$ and $P_n(\Lambda)$. This is the main result of this section.
\begin{thm}\label{cor03} Let $A$ be a finite-dimensional Hopf algebra over the field $\kk$ with  a nonzero left integral $\Lambda$ and a normalized twist $J$.  For any $n\in\mathbb{Z}$, we have
\begin{enumerate}
  \item $P_n^J(\Lambda)=TP_n(\Lambda)$, where $T=J^{-(1)}\alpha(J^{-(2)})$.
  \item $P_n^J(\Lambda)=Q_J^{-1}P_{n}(\Lambda)S(R)Q_J$, where $R=\alpha(J^{-(1)})J^{-(2)}$.
\end{enumerate}
If, moreover, $A$ is unimodular, then $P_n^J(\Lambda)=P_n(\Lambda)$.
\end{thm}
\proof (1)
It is direct to check that the identity in question holds for the case $n=0,\pm1$. For the case $n\geq2$, we have
\begin{align*}
P_n^J(\Lambda)&=\Lambda_{\langle1\rangle}\Lambda_{\langle2\rangle}\Lambda_{\langle3\rangle}\cdots\Lambda_{\langle n\rangle}\\
&=Q_J^{-1}\Lambda_{(1)}J^{(1)}J^{-(1)}\Lambda_{(2)}P_{n-2}\bigg(J^{(2)}J^{-(2)}\Lambda_{(3)}\bigg)S(R)Q_J\ \ \ \text{by\ Lemma}\ \ref{prop1}\ (1)\\
&=Q_J^{-1}\Lambda_{(1)}\Lambda_{(2)}P_{n-2}(\Lambda_{(3)})S(R)Q_J\\
&=Q_J^{-1}P_{n}(\Lambda)S(R)Q_J\\
&=P_{n}(\Lambda)J^{(1)}(\alpha\circ S^{-1})(J^{(2)}Q_J^{-1})\ \ \ \text{by\ Proposition}\ \ref{lem1}\\
&=P_{n}(\Lambda)J^{-(1)}_{(1)}(\alpha\circ S^{-1})(J^{-(1)}_{(2)})\alpha(J^{-(2)})\ \ \ \text{by}\ (\ref{equ0033})\\
&=TP_n(\Lambda)\ \ \ \text{by\ Proposition}\ \ref{lem1}
\end{align*}
For the case $n\leq -2$, we have
\begin{align*}
P_n^J(\Lambda)&=S^J(\Lambda_{\langle-n\rangle}\cdots\Lambda_{\langle1\rangle})\\
&=Q_J^{-1}S\bigg(\Lambda_{\langle-n\rangle}\cdots\Lambda_{\langle1\rangle}\bigg)Q_J\\
&=Q_J^{-1}S\bigg(S^{-1}(Q_J^{-1})\Lambda_{(3)}J^{(2)}J^{-(2)}\Lambda_{(2)}(S^{-1}\circ P_{n+2})(J^{(1)}J^{-(1)}\Lambda_{(1)})S^{-1}(T)S^{-1}(Q_J)\bigg)Q_J\\
&=Q_J^{-1}S\bigg(S^{-1}(Q_J^{-1})\Lambda_{(3)}\Lambda_{(2)}(S^{-1}\circ P_{n+2})(\Lambda_{(1)})S^{-1}(T)S^{-1}(Q_J)\bigg)Q_J\\
&=S\bigg(\Lambda_{(3)}\Lambda_{(2)}(S^{-1}\circ P_{n+2})(\Lambda_{(1)})S^{-1}(T)\bigg)\\
&=TP_n(\Lambda),
\end{align*}
where the third equality follows from  Lemma \ref{prop1} (2).

(2) On the one hand,
\begin{align*}
P_n^J(\Lambda)&=TP_n(\Lambda)\\
&=P_n(\Lambda)J^{-(1)}_{(1)}(\alpha\circ S^{-1})(J^{-(1)}_{(2)})\alpha(J^{-(2)})\ \ \ \text{by\ Proposition}\ \ref{lem1}\\
&=P_{n}(\Lambda)J^{(1)}(\alpha\circ S^{-1})(J^{(2)}Q_J^{-1}) \ \ \text{by}\ (\ref{equ0033})
\end{align*}
On the other hand,
\begin{align*}
Q_J^{-1}P_{n}(\Lambda)S(R)Q_J&=P_{n}(\Lambda)(Q_J^{-1})_{(1)}(\alpha\circ S^{-1})((Q_J^{-1})_{(2)})S(R)Q_J \ \ \text{by\ Proposition}\ \ref{lem1}\\
&=P_n(\Lambda)J^{(1)}Q_J^{-1}S(\widetilde{J}^{(2)})(\alpha\circ S^{-1})(J^{(2)}Q_J^{-1}S(\widetilde{J}^{(1)}))S(R)Q_J\ \text{by}\ (\ref{equu00003})\\
&=P_{n}(\Lambda)J^{(1)}(\alpha\circ S^{-1})(J^{(2)}Q_J^{-1}).
\end{align*}We conclude that $P_n^J(\Lambda)=Q_J^{-1}P_{n}(\Lambda)S(R)Q_J$ for any $n\in\mathbb{Z}$.

If $A$ is unimodular, then $\alpha=\varepsilon$ and hence $T=1$, this implies that $P_n^J(\Lambda)=P_n(\Lambda)$ for any $n\in\mathbb{Z}$.
\qed

\begin{rem}
Theorem \ref{cor03} shows that the sequence $\{P_n(\Lambda)\}_{n\in\mathbb{Z}}$ of a unimodular Hopf algebra is invariant under twisting.
\end{rem}


\section{\bf Invariants of representation categories of Hopf algebras}
In this section, we will use Theorem \ref{cor03} to give several invariants of representation categories of Hopf algebras. As applications, we show that the representation categories of 12-dimensional pointed nonsemisimple Hopf algebras classified in \cite{Natale} are mutually inequivalent as tensor categories. We also show that those $8$-dimensional semisimple Hopf algebras, including Kac algebra $K_8$,  the dihedral group algebra $\kk D_4$ and the  quaternion group algebra $\kk Q_8$ are mutually not gauge equivalent.

We begin with the following preparations. Let $A$ and $A'$ be finite-dimensional Hopf algebras over the field $\kk$ with nonzero left integrals $\Lambda$ and $\Lambda'$ respectively. If the functor $\mathcal{F}:\text{Rep}(A)\rightarrow\text{Rep}(A')$ is an equivalence of tensor categories, it follows from \cite[Theorem 2.2]{NS1} that there exist
a normalized twist $J$ of $A$  such that $A'$ is isomorphic to $A^J$ as bialgebras. Let $\sigma:A'\rightarrow A^J$ be such an isomorphism. Then $\sigma$ is
automatically a Hopf algebra isomorphism. Therefore, for any $n\in\mathbb{Z}$, we have
$$\sigma\circ P'_n=P_n^J\circ\sigma,$$ where $P'_n$ and $P_n^J$ are the $n$-th Sweedler power maps of $A'$ and $A^J$ respectively.
Suppose $\sigma(\Lambda')=\mu\Lambda$ for a nonzero scalar $\mu\in\kk$. Then
\begin{equation}\label{qu01}\sigma(P'_n(\Lambda'))
=P^J_n(\sigma(\Lambda'))=\mu P^J_n(\Lambda)=\mu TP_n(\Lambda),
\end{equation}
where $T=J^{-(1)}\alpha(J^{-(2)})$.
The isomorphism $\sigma$ induces a $\kk$-linear equivalence $(-)^\sigma:\text{Rep}(A)\rightarrow\text{Rep}(A')$ as follows:  for any finite-dimensional $A$-module $V$, $V^{\sigma}=V$ as $\kk$-linear space with the $A'$-module structure given by $a'v=\sigma(a')v$ for $a'\in A'$, $v\in V$, and $f^\sigma=f$ for any morphism $f$ in Rep$(A)$.
Thus, $$\chi_{V^{\sigma}}(a')=\chi_V(\sigma(a'))\ \text{for}\ a'\in A'.$$
Moreover, the equivalence  $\mathcal{F}$ is naturally isomorphic to the $\kk$-linear equivalence $(-)^\sigma$ (see \cite[Theorem 1.1]{KMN}).
Therefore, $$\chi_{\mathcal{F}(V)}(a')=\chi_{V^{\sigma}}(a')\ \text{for}\ a'\in A'.$$
By taking $a'$ to be $P'_n(\Lambda')$, we have
\begin{equation}\label{qqu}\chi_{\mathcal{F}(V)}(P'_n(\Lambda'))=\chi_{V^{\sigma}}(P'_n(\Lambda'))=\chi_{V}(\sigma(P'_n(\Lambda')))=\mu\chi_{V}(TP_n(\Lambda)),\end{equation}
where the scalar $\mu$ is independent on the choice of $V$.

The first result of this section states that any homogeneous polynomial vanishing on some $P_n(\Lambda)$ is an invariant of the tensor category Rep$(A)$. We shall see that this invariant can be used to distinguish the representation categories of some Hopf algebras.

\begin{thm}\label{p1}
Let $A$ and $A'$ be gauge equivalent finite-dimensional Hopf
algebras over the field $\kk$ with nonzero left integrals $\Lambda$ and $\Lambda'$ respectively. For any homogeneous polynomial $\psi(X_1,\cdots,X_s)\in\kk[X_1,\cdots,X_s]$, we have  $\psi(P_{n_1}(\Lambda),\cdots,P_{n_s}(\Lambda))=0$ for some $n_1,\cdots,n_s\in\mathbb{Z}$ if and only if $\psi(P'_{n_1}(\Lambda'),\cdots,P'_{n_s}(\Lambda'))=0$.
\end{thm}
\proof Since Rep$(A)$ and Rep$(A')$ are equivalent as tensor categories, there exists
a normalized twist $J$ of $A$ such that $A'$ is isomorphic to $A^J$ as Hopf algebras. Let $\sigma:A'\rightarrow A^J$ be such an isomorphism. It follows from (\ref{qu01}) that $\sigma(P'_n(\Lambda'))
=\mu TP_n(\Lambda)$, where $\mu$ is a nonzero scalar and $T=J^{-(1)}\alpha(J^{-(2)})$ which is an invertible element of $A$. Let $\psi(X_1,\cdots,X_s)$ be a homogeneous polynomial of degree $m$. Then
\begin{align*}
\sigma(\psi(P'_{n_1}(\Lambda'),\cdots,P'_{n_s}(\Lambda')))&=\psi(\sigma(P'_{n_1}(\Lambda')),\cdots,\sigma(P'_{n_s}(\Lambda')))\\
&=\psi(\mu TP_{n_1}(\Lambda),\cdots,\mu TP_{n_s}(\Lambda))\\
&=\psi(P_{n_1}(\Lambda),\cdots,P_{n_s}(\Lambda))\mu^{m}T_{m}\cdots T_2T_1\ \ \text{by}\ (\ref{eer})
\end{align*}
where these elements $T_1,\cdots,T_m\in A$ are as defined recursively in (\ref{er}). Note that these elements $T_1,\cdots,T_m$ are invertible in $A$.
Therefore, $\psi(P_{n_1}(\Lambda),\cdots,P_{n_s}(\Lambda))=0$ if and only if $\psi(P'_{n_1}(\Lambda'),\cdots,P'_{n_s}(\Lambda'))=0$. The proof is completed.
\qed

We use Theorem \ref{p1} to distinguish the representation categories of $12$-dimensional pointed nonsemisimple Hopf algebras. The following is a list of pairwise nonisomorphic pointed nonsemisimple Hopf algebras of
dimension 12 over an algebraically closed field $\kk$ of characteristic zero. Every Hopf algebra in the list is presented by two generators $g$ and $x$ subject to the following relations:
\begin{align*}
&\mathcal{A}_0:\ \ g^6=1,\ x^2=0,\ gx=-xg,\ \Delta(g)=g\otimes g,\ \Delta(x)=x\otimes1+g\otimes x,\\
&\mathcal{A}_1:\ \ g^6=1,\ x^2=1-g^2,\ gx=-xg,\ \Delta(g)=g\otimes g,\ \Delta(x)=x\otimes1+g\otimes x,\\
&\mathcal{B}_0:\ \ g^6=1,\ x^2=0,\ gx=-xg,\ \Delta(g)=g\otimes g,\ \Delta(x)=x\otimes1+g^3\otimes x,\\
&\mathcal{B}_1:\ \ g^6=1,\ x^2=0,\ gx=\omega xg,\ \Delta(g)=g\otimes g,\ \Delta(x)=x\otimes1+g^3\otimes x,
\end{align*}
where $\omega\in\kk$ is a fixed primitive $6$-th root of unity. It follows from \cite{Natale} that these four Hopf algebras are up to isomorphism all the pointed nonsemisimple Hopf algebras of dimension 12. These Hopf algebras are not unimodular and have nonzero left integrals of the same form
$$\Lambda=(1+g+g^2+g^3+g^4+g^5)x.$$ Since $S^2(\Lambda)=-\Lambda$ holds for all these Hopf algebras, it follows from \cite[Proposition 3.13]{Sh} that the $(-1)$-th indicators $\nu_{-1}$ of these Hopf algebras are all $-1$. For the trace of antipode $S$ of each Hopf algebra, we have tr$(S)=2$. Hence the gauge invariants $\nu_{-1}$ and $\nu_2=\text{tr}(S)$ which are easy to handle can not be used to distinguish the representation categories of these Hopf algebras. However, we may find some homogeneous polynomials of degree 1 vanishing on some $P_n(\Lambda)$ as follows:
\begin{table}[H]
\begin{tabular}{|c|c|c|c|}
\hline
$\mathcal{A}_0$ & $P_2(\Lambda)+P_{-2}(\Lambda)=0$ & $P_3(\Lambda)-3P_{2}(\Lambda)-P_{-3}(\Lambda)=0$  \\
\hline
$\mathcal{A}_1$ & $P_2(\Lambda)+P_{-2}(\Lambda)=0$ & $P_3(\Lambda)-3P_{2}(\Lambda)-P_{-3}(\Lambda)=0$  \\
\hline
$\mathcal{B}_0$ & $P_2(\Lambda)+P_{-2}(\Lambda)=0$ & $P_3(\Lambda)-3P_{2}(\Lambda)-P_{-3}(\Lambda)\neq0$   \\
\hline
$\mathcal{B}_1$ & $P_2(\Lambda)+P_{-2}(\Lambda)\neq0$ &   \\
\hline
\end{tabular}
\end{table}
According to Theorem \ref{p1}, we may see from the second column of the above table that $\mathcal{B}_1$ is not gauge equivalent to  $\mathcal{A}_0$, $\mathcal{A}_1$ and $\mathcal{B}_0$. Similarly, it can be seen from the third column that  $\mathcal{B}_0$ is not gauge equivalent to  $\mathcal{A}_0$ and $\mathcal{A}_1$. However, this approach can not be used to distinguish the representation categories $\mathcal{A}_0$ and $\mathcal{A}_1$, since $P_n(\Lambda)\in\mathcal{A}_0$ and $P_n(\Lambda)\in\mathcal{A}_1$ have the same expression for any $n\in\mathbb{Z}$ with respect to the $\kk$-basis $\{g^ix^j\mid 0\leq i\leq5,0\leq j\leq1\}$ of $\mathcal{A}_0$ and of $\mathcal{A}_1$. Fortunately, the representation categories of $\mathcal{A}_0$ and $\mathcal{A}_1$ have already been investigated in\cite{WLZ1} and \cite{WLZ2} respectively. The number of finite-dimensional indecomposable representations of $\mathcal{A}_0$ up to isomorphism is 12 (see \cite[Theorem 2.5]{WLZ1}), while the number for that of $\mathcal{A}_1$ is 6 (see \cite[Theorem 2.9]{WLZ2}).
In summary, the Hopf algebras $\mathcal{A}_0$, $\mathcal{A}_1$, $\mathcal{B}_0$ and $\mathcal{B}_1$ are mutually not gauge equivalent.

Note that the (not necessarily homogeneous) polynomials vanishing on some $P_n(\Lambda)$ are in general not anymore invariants of the tensor category Rep$(A)$. But for a unimodular Hopf algebra $A$, as shown in the following that these polynomials are still something meaningful to distinguish the tensor category Rep$(A)$.

\begin{prop}\label{p2}Let $A$ and $A'$ be gauge equivalent finite-dimensional unimodular Hopf
algebras over the field $\kk$. Let  $\Lambda$  be a nonzero left integral in $A$.
For any polynomial $\psi(X_1,\cdots,X_s)\in\kk[X_1,\cdots,X_s]$, if $\psi(P_{n_1}(\Lambda),\cdots,P_{n_s}(\Lambda))=0$ for some $n_1,\cdots,n_s\in\mathbb{Z}$, then there exists a nonzero left integral $\Lambda'$ in $A'$ such that $\psi(P'_{n_1}(\Lambda'),\cdots,P'_{n_s}(\Lambda'))=0$.
\end{prop}
\proof
Since Rep$(A)$ and Rep$(A')$ are equivalent as tensor categories, there exists
a normalized twist $J$ of $A$ such that $A'$ is isomorphic via a map $\sigma$ to $A^J$ as Hopf algebras. We may choose a left integral $\Lambda'\in A'$ such that $\sigma(\Lambda')=\Lambda$. In this case, the scalar $\mu$ appeared in (\ref{qu01}) is 1. Since $A$ is unimodular, the element $T$ appeared in (\ref{qu01}) is also 1. Now (\ref{qu01}) has the form
$\sigma(P'_n(\Lambda'))=P_n(\Lambda)$ for any $n\in\mathbb{Z}$. Thus, $\psi(P_{n_1}(\Lambda),\cdots,P_{n_s}(\Lambda))=0$ implies that $\psi(P'_{n_1}(\Lambda'),\cdots, P'_{n_s}(\Lambda'))=0$. The proof is completed.
\qed

For a semisimple Hopf algebra $A$ with an idempotent integral $\Lambda$, as shown in the following that any polynomial vanishing on some $P_n(\Lambda)$ is indeed an invariant of the tensor category Rep$(A)$.

\begin{prop}\label{r1}Let $A$ and $A'$ be gauge equivalent semisimple Hopf algebras with idempotent integrals $\Lambda$ and $\Lambda'$ respectively. For any polynomial $\psi(X_1,\cdots,X_s)\in\kk[X_1,\cdots,X_s]$, we have  $\psi(P_{n_1}(\Lambda),\cdots,P_{n_s}(\Lambda))=0$ for some $n_1,\cdots,n_s\in\mathbb{Z}$ if and only if $\psi( P'_{n_1}(\Lambda'),\cdots,P'_{n_s}(\Lambda'))=0$.
\end{prop}
\proof Since Rep$(A)$ and Rep$(A')$ are equivalent as tensor categories, there exists
a normalized twist $J$ of $A$ such that $A'$ is isomorphic via a map $\sigma$ to $A^J$ as Hopf algebras. Since $\Lambda$ and $\Lambda'$ are both idempotent and $A$ is unimodular, the scalar $\mu$ and the element $T$ appeared in (\ref{qu01}) are both equal to $1$. Now (\ref{qu01}) becomes $\sigma(P'_n(\Lambda'))=P_n(\Lambda)$ for any $n\in\mathbb{Z}$. Thus, $\psi(P_{n_1}(\Lambda),\cdots,P_{n_s}(\Lambda))=0$ if and only if $\psi(P'_{n_1}(\Lambda'),\cdots,P'_{n_s}(\Lambda'))=0$. The proof is completed.
\qed

In the sequel, we will use Proposition \ref{r1} to distinguish the representation categories Rep$(\kk Q_8)$, Rep$(\kk D_4)$ and Rep$(K_8)$ over the field $\kk$ of characteristic $\neq2$.
Note that these representation categories have the same Grothendieck ring but they are inequivalent as tensor categories as they have different higher Frobenius-Schur indicators (see \cite[Theorem 6.1]{NS1}).

The quaternion group $Q_8$ is a group with eight elements, which can be described in the following way:
It is the group formed by eight elements $\mathbbm{1},-\mathbbm{1},i,-i,j,-j,k,-k$ where $\mathbbm{1}$ is the identity element, $(-\mathbbm{1})^2=\mathbbm{1}$ and all the other elements are squareroots of $-\mathbbm{1}$, such that $(-\mathbbm{1})i=-i, (-\mathbbm{1})j=-j, (-\mathbbm{1})k=-k$ and further, $ij=k, ji=-k, jk=i, kj=-i, ki=j, ik=-j$ (the remaining relations can be deduced from these). The group algebra $\kk Q_8$ has the idempotent integral $$\Lambda=\frac{1}{8}(\mathbbm{1}+(-\mathbbm{1})+i+(-i)+j+(-j)+k+(-k)).$$ By a straightforward computation, we have
$$P_1(\Lambda)=\Lambda,\ \ \ \ P_2(\Lambda)=\frac{1}{4}\mathbbm{1}+\frac{3}{4}(-\mathbbm{1}),\ \ \  P_3(\Lambda)=P_1(\Lambda),\ \ \ \ P_4(\Lambda)=P_0(\Lambda).$$
Note that $P_4(\Lambda)-P_0(\Lambda)=0$ and $2(P_2(\Lambda))^2-P_2(\Lambda)-P_0(\Lambda)=0$.

For the dihedral group $D_4 =\{ 1,a,a^2,a^{3},b,ba,ba^2,ba^{3}\}$ with $a^4=1,\ b^2=1$ and $aba=b$, the group algebra  $\kk D_4$ has the idempotent integral $$\Lambda=\frac{1}{8}(1+a+a^2+a^3+b+ba+ba^2+ba^3).$$ We have
$$P_1(\Lambda)=\Lambda,\ \ \ \ P_2(\Lambda)=\frac{3}{4}+\frac{1}{4}a^2,\ \ \ \ P_3(\Lambda)=P_1(\Lambda),\ \ \ \ P_4(\Lambda)=P_0(\Lambda).$$
Note that $P_4(\Lambda)-P_0(\Lambda)=0$ while $2(P_2(\Lambda))^2-P_2(\Lambda)-P_0(\Lambda)=\frac{1}{2}a^2-\frac{1}{2}\neq0$

The 8-dimensional Kac algebra $K_8$ is a semisimple Hopf algebra over $\kk$ generated
by $x, y, z$ as a $\kk$-algebra with the following relations (see \cite{Mas}):
$$x^2=y^2=1,\ z^2=\frac{1}{2}(1+x+y-xy),\ xy=yx,\ xz=zy,\ yz=zx.$$
The coalgebra structure $\Delta,\varepsilon$ and the antipode $S$ of $K_8$ are given by
$$\Delta(x)=x\otimes x,\ \Delta(y)=y\otimes y,\ \varepsilon(x)=\varepsilon(y)=1,$$
$$\Delta(z)=\frac{1}{2}(1\otimes1+1\otimes x+y\otimes1-y\otimes x)(z\otimes z),\ \varepsilon(z)=1,$$
$$S(x)=x,\ S(y)=y,\ S(z)=z.$$
The idempotent integral of $K_8$ is $$\Lambda=\frac{1}{8}(1+x+y+xy)+\frac{1}{8}(1+x+y+xy)z.$$ A straightforward computation shows that
$$P_1(\Lambda)=\Lambda,\ \ \ \ P_2(\Lambda)=\frac{3}{4}+\frac{1}{4}xy,\ \ \ \ P_3(\Lambda)=P_1(\Lambda),\ \ \ \ P_4(\Lambda)=\frac{1}{2}+\frac{1}{2}xy,$$
$$P_5(\Lambda)=P_1(\Lambda),\ \ \ \ P_6(\Lambda)=P_2(\Lambda),\ \ \ \ P_7(\Lambda)=P_1(\Lambda),\ \ \ \ P_8(\Lambda)=P_0(\Lambda).$$
Note that $P_4(\Lambda)-P_0(\Lambda)=\frac{1}{2}xy-\frac{1}{2}\neq0$.

We may see from the following table that for any two of these  Hopf algebras, there exists a polynomial $\psi$ in three variables such that $\psi(P_0(\Lambda),P_2(\Lambda),P_4(\Lambda))=0$ holds for one Hopf algebra but not for another Hopf algebra. Thus, the three representation categories are mutually inequivalent as tensor categories by Proposition \ref{r1}.

\begin{table}[H]
\begin{tabular}{|c|c|c|c|}
\hline
$K_8$ & $P_4(\Lambda)-P_0(\Lambda)\neq0$ &  & $P_4(\Lambda)-P_0(\Lambda)\neq0$ \\
\hline
$\kk D_4$ & $P_4(\Lambda)-P_0(\Lambda)=0$ & $2(P_2(\Lambda))^2-P_2(\Lambda)-P_0(\Lambda)\neq0$ &   \\
\hline
$\kk Q_8$ &  & $2(P_2(\Lambda))^2-P_2(\Lambda)-P_0(\Lambda)=0$ & $P_4(\Lambda)-P_0(\Lambda)=0$ \\
\hline
\end{tabular}
\end{table}


Although these $P_n(\Lambda)$ for $n\in\mathbb{Z}$ are invariant under twisting for a finite-dimensional unimodular Hopf algebra $A$, it is clear that  $P_n(\Lambda)$ are dependent on the choice of $\Lambda$, and the values that characters of finite-dimensional $A$-modules taking on them are not gauge invariants in general. Note that if $A$ is unimodular then the element $T$ appeared in (\ref{qqu}) is the identity 1. Now (\ref{qqu}) has the form $$\chi_{\mathcal{F}(V)}(P'_n(\Lambda'))=\mu\chi_{V}(P_n(\Lambda))\ \text{for}\ n\in\mathbb{Z},$$ where $\mu$ is a nonzero scalar which is independent on the choice of $V$. In view of this, if $\chi_{W}(P_m(\Lambda))\neq0$ for some finite-dimensional $A$-module $W$ and $m\in\mathbb{Z}$, then for any finite-dimensional $A$-module $V$ and $n\in\mathbb{Z}$, we have
$$\chi_{\mathcal{F}(V)}(P'_n(\Lambda')):\chi_{\mathcal{F}(W)}(P'_m(\Lambda'))=\chi_{V}(P_n(\Lambda)):\chi_{W}(P_m(\Lambda)).$$
That is, the ratio
$\chi_{V}(P_n(\Lambda)):\chi_{W}(P_m(\Lambda))$ is a gauge invariant of a finite-dimensional unimodular Hopf algebra $A$. We summarize it as follows:

\begin{prop}\label{prop112}
Let $A$ be a finite-dimensional unimodular Hopf algebra over the field $\kk$ with a nonzero left integral $\Lambda$.
If $\chi_{W}(P_m(\Lambda))\neq0$ for some finite-dimensional $A$-module $W$ and $m\in\mathbb{Z}$, then for any finite-dimensional $A$-module $V$ and $n\in\mathbb{Z}$,  the ratio
$\chi_{V}(P_n(\Lambda)):\chi_{W}(P_m(\Lambda))$ is a gauge invariant of $A$.
\end{prop}

For a semisimple Hopf algebra $A$ over the field $\kk$ with an idempotent integral $\Lambda$, the value $\chi_V(P_n(\Lambda))$ is called the $n$-th FS-indicator of a finite-dimensional $A$-module $V$ (see \cite{KSZ,LM}). It follows from \cite[Proposition 3.2]{NS1} that $\chi_V(P_n(\Lambda))$  is a gauge invariant of $A$ over the field of complex numbers. Proposition \ref{prop112} can be regarded as a slightly generalization of the $n$-th FS-indicator from a semisimple Hopf algebra to a unimodular Hopf algebra. Indeed, if $A$ is semisimple and $W$ is chosen to be the trivial $A$-module $\kk$, then $\chi_{\kk}(P_m(\Lambda))=\varepsilon(\Lambda)\neq0$ and the ratio $\chi_{V}(P_n(\Lambda)):\chi_{\kk}(P_m(\Lambda))=\chi_{V}(P_n(\Lambda/\varepsilon(\Lambda)))$ is the $n$-th FS-indicator of $V$ for any $n\in\mathbb{Z}$. Note that the base field $\kk$ is arbitrary. It suggests that the invariant of the $n$-th FS-indicator of $V$ is valid for a semisimple Hopf algebra over an arbitrary field $\kk$.


For a finite-dimensional non-semisimple Hopf algebra $A$, the following result suggests that it is interesting to seek for an $A$-module $W$ such that $\chi_{W}(P_m(\Lambda))\neq0$ for some $m\in\mathbb{Z}$.

\begin{prop}
Let $A$ be a finite-dimensional non-semisimple Hopf algebra with the Chevalley property. For any finite-dimensional $A$-module $V$ and $n\in\mathbb{Z}$, we always have $\chi_V(P_n(\Lambda))=0$.
\end{prop}
\proof
Since $A$ is a non-semisimple Hopf algebra, any nonzero left integral $\Lambda$ is a nilpotent element of $A$. It follows that $\Lambda$ belongs to the Jacobson radical $\text{rad}(A)$ of $A$. Note that the Hopf algebra $A$ has the Chevalley property. Namely, the Jacobson radical $\text{rad}(A)$ is a Hopf ideal of $A$. Therefore, $P_n(\Lambda)\in\text{rad}(A)$ and hence $\chi_V(P_n(\Lambda))=0$ for any $n\in\mathbb{Z}$. \qed

At the end of this section, we pay a little attention to the Hopf order of a nonzero left integral of $A$. Recall from \cite{LMS} that the Hopf order of a nonzero left integral $\Lambda\in A$ is the least positive integer $n$ such that $P_n(\Lambda)=P_0(\Lambda)$. If such $n$ does not exist, the Hopf order of $\Lambda$ is infinity.
If $A$ and $A'$ are gauge equivalent Hopf algebras with  nonzero left integrals $\Lambda$ and $\Lambda'$ respectively, then $\sigma(P'_n(\Lambda'))
=\mu TP_n(\Lambda)$ by (\ref{qu01}). It follows that
 $P'_n(\Lambda')=P'_0(\Lambda')$ if and only if $P_n(\Lambda)=P_0(\Lambda)$. Therefore, the Hopf order of $\Lambda$ is equal to that of $\Lambda'$. That is, the Hopf order of a nonzero left integral $\Lambda\in A$ is a gauge invariant of $A$.
It can be seen from above that the Hopf orders of $\Lambda\in \kk D_4$ and $\Lambda\in \kk Q_8$ are both equal to 4, while the Hopf order of $\Lambda\in K_8$ is equal to 8.


\section{\bf Integrals in the dual of twisted Hopf algebras}
In this section, we investigate the relationship between a right integral $\lambda\in A^*$ and a right integral $\lambda^J\in(A^J)^*$. Based on this investigation, we provide a unifying proof of the well-known result which says that the indicator $\nu_n(A)$ of a finite-dimensional Hopf algebra $A$ is a gauge invariant of $A$ for any $n\in\mathbb{Z}$. We also use the expression of $\lambda^J$ to give a different proof of the known result that the Killing form of a finite-dimensional Hopf algebra $A$ is invariant under twisting.

Note that the expression of a right integral $\lambda^J\in(A^J)^*$ has been described in \cite[Theorem 3.4]{AEGN} when $A$ is a finite-dimensional unimodular Hopf algebra. For a general Hopf algebra $A$, a right integral $\lambda^J\in(A^J)^*$ can be described as follows:

\begin{thm}\label{th2}
Let $A$ be a finite-dimensional Hopf algebra over the field $\kk$ with a normalized twist $J$. Let $R=\alpha(J^{-(1)})J^{-(2)}$. If $\lambda$ is a nonzero right integral in $A^*$, then $$\lambda^J:=\lambda\leftharpoonup S^2(R^{-1})S(Q_J^{-1})Q_J$$ is a nonzero right integral in $(A^J)^*$.
\end{thm}
\proof We need to show that $\lambda^J(b_{\langle1\rangle})b_{\langle2\rangle}=\lambda^J(b)$ for all $b\in A$. We set $\widetilde{J}^{(1)}\otimes \widetilde{J}^{(2)}=J^{(1)}\otimes J^{(2)}$. Then
\begin{align*}\lambda^J(b_{\langle1\rangle})b_{\langle2\rangle}&=\lambda^J\bigg(J^{-(1)}b_{(1)}J^{(1)}\bigg)J^{-(2)}b_{(2)}J^{(2)}\\
&=\lambda\bigg(\alpha(\widetilde{J}^{(1)})S^2(\widetilde{J}^{(2)})S(Q_J^{-1})Q_JJ^{-(1)}b_{(1)}J^{(1)}\bigg)J^{-(2)}b_{(2)}J^{(2)}\\
&=\lambda\bigg(\alpha(J^{(1)}_{(1)})S^2(J^{(1)}_{(2)})\alpha(\widetilde{J}^{(1)})S^2(\widetilde{J}^{(2)})S(Q_J^{-1})Q_JJ^{-(1)}b_{(1)}\bigg)J^{-(2)}b_{(2)}J^{(2)}\ \ \text{by}\ (\ref{equ4004})\\
&=\lambda\bigg(\alpha(J^{(1)}_{(1)}\widetilde{J}^{(1)})S^2(J^{(1)}_{(2)}\widetilde{J}^{(2)})S(Q_J^{-1})Q_JJ^{-(1)}b_{(1)}\bigg)J^{-(2)}b_{(2)}J^{(2)}\\
&=\lambda\bigg(\alpha(J^{(1)})S^2(J^{(2)}_{(1)}\widetilde{J}^{(1)})S(Q_J^{-1})Q_JJ^{-(1)}b_{(1)}\bigg)J^{-(2)}b_{(2)}J^{(2)}_{(2)}\widetilde{J}^{(2)}\ \ \text{by}\ (\ref{qq1})
\end{align*}
We denote $t:=S(Q_J^{-1})Q_J.$ It follows from (\ref{equ002}) and (\ref{equu00003}) that
$$\Delta(t)=t_{(1)}\otimes t_{(2)}=(S^2\otimes S^2)(J)(S(Q_J^{-1})Q_J\otimes S(Q_J^{-1})Q_J)J^{-1}.$$ Applying $id\otimes S^{-1}$ to both sides of this equality, we obtain that
\begin{equation}\label{r}t_{(1)}\otimes S^{-1}(t_{(2)})
=S^2(\widetilde{J}^{(1)})S(Q_J^{-1})Q_J J^{-(1)}\otimes S^{-1}(J^{-(2)})S^{-1}(Q_J)Q_J^{-1}S(\widetilde{J}^{(2)}).
\end{equation}
Now we have
\begin{align*}&\ \lambda^J(b_{\langle1\rangle})b_{\langle2\rangle}\\
&=\alpha(J^{(1)})J^{-(2)}\lambda\bigg(S^2(J^{(2)}_{(1)}\widetilde{J}^{(1)})tJ^{-(1)}b_{(1)}\bigg)b_{(2)}J^{(2)}_{(2)}\widetilde{J}^{(2)}\\
&=\alpha(J^{(1)})J^{-(2)}\lambda\bigg(S^2(J^{(2)}_{(1)}\widetilde{J}^{(1)}_{(1)})t_{(1)}J^{-(1)}_{(1)}b\bigg)
S^{-1}\bigg(S^2(J^{(2)}_{(2)}\widetilde{J}^{(1)}_{(2)})t_{(2)}J^{-(1)}_{(2)}\bigg)J^{(2)}_{(3)}\widetilde{J}^{(2)}\ \ \text{by}\ (\ref{equ4004444})\\
&=\alpha(J^{(1)})J^{-(2)}\lambda\bigg(S^2(J^{(2)}_{(1)}\widetilde{J}^{(1)}_{(1)})t_{(1)}J^{-(1)}_{(1)}b\bigg)
S^{-1}(J^{-(1)}_{(2)})S^{-1}(t_{(2)})S(\widetilde{J}^{(1)}_{(2)})S(J^{(2)}_{(2)})
J^{(2)}_{(3)}\widetilde{J}^{(2)}\\
&=\alpha(J^{(1)})J^{-(2)}\lambda\bigg(S^2(J^{(2)}\widetilde{J}^{(1)}_{(1)})t_{(1)}J^{-(1)}_{(1)}b\bigg)
 S^{-1}(J^{-(1)}_{(2)})S^{-1}(t_{(2)})S(\widetilde{J}^{(1)}_{(2)})
\widetilde{J}^{(2)}\\
&=\lambda\bigg(S^2(R^{-1})S^2(\widetilde{J}^{(1)}_{(1)})t_{(1)}J^{-(1)}_{(1)}b\bigg)
J^{-(2)}S^{-1}(J^{-(1)}_{(2)})S^{-1}(t_{(2)})S(\widetilde{J}^{(1)}_{(2)})
\widetilde{J}^{(2)}\\
&=\lambda\bigg(S^2(R^{-1})S^2(\widetilde{J}^{(1)}_{(1)})t_{(1)}J^{(1)}b\bigg)
S^{-1}(Q_J^{-1})S^{-1}(J^{(2)})S^{-1}(t_{(2)})S(\widetilde{J}^{(1)}_{(2)})
\widetilde{J}^{(2)}\ \ \text{by}\ (\ref{equ0033})\\
&=\lambda\bigg(S^2(R^{-1})S^2(\widetilde{J}^{-(1)})t_{(1)}J^{(1)}b\bigg)
S^{-1}(Q_J^{-1})S^{-1}(J^{(2)})S^{-1}(t_{(2)})S(\widetilde{J}^{-(2)})Q_J\ \text{by}\ (\ref{equ003})\\
&=\lambda\bigg(S^2(R^{-1})S^2(\widetilde{J}^{-(1)})S^2(\widetilde{J}^{(1)})S(Q_J^{-1})Q_J J^{-(1)}J^{(1)}b\bigg)\\
&\ \ \ \ \ \ \ \ \ \ \ \ \ \cdot S^{-1}(Q_J^{-1})S^{-1}(J^{(2)})S^{-1}(J^{-(2)})S^{-1}(Q_J)Q_J^{-1}S(\widetilde{J}^{(2)})S(\widetilde{J}^{-(2)})Q_J\ \ \text{by}\ ({\ref{r}})\\
&=\lambda\bigg(S^2(R^{-1})S(Q_J^{-1})Q_J b\bigg)S^{-1}(Q_J^{-1})S^{-1}(Q_J)Q_J^{-1}Q_J\\
&=\lambda\bigg(S^2(R^{-1})S(Q_J^{-1})Q_J b\bigg)=\lambda^J(b).
\end{align*}
Thus, $\lambda^J$ is a right integral in $(A^J)^*$. The proof is completed.
\qed

\begin{rem}
If $A$ is a finite-dimensional unimodular Hopf algebra over the field $\kk$ with a normalized twist $J$, then $R=\alpha(J^{-(1)})J^{-(2)}=\varepsilon(J^{-(1)})J^{-(2)}=1$. In this case, $\lambda^J=\lambda\leftharpoonup S(Q_J^{-1})Q_J$, which is exactly the result of \cite[Theorem 3.4]{AEGN}. It is interesting to know whether or not $\lambda(1)\neq0$ implying $\lambda^J(1)\neq0$ for a finite-dimensional Hopf algebra $A$ (see \cite[Remark 3.9]{AEGN}). This problem has only been solved in the unimodular case (see \cite[Corollary 3.6]{AEGN}).
\end{rem}

For any finite-dimensional Hopf algebra $A$, recall from \cite{KMN} and \cite{Sh} that the $n$-th indicator of $A$ is defined by
 $$\nu_n(A):=\text{tr}(S\circ P_{n-1})\ \ \text{for}\ n\in\mathbb{Z}.$$

By Radford's trace formula (\ref{equ400444}), we have
\begin{equation}\label{equ01}\nu_n(A)=\left\{
             \begin{array}{ll}
               (\lambda\circ S)(\Lambda_{(1)}\Lambda_{(2)}\cdots\Lambda_{(n)}), & n\geq1; \\
               \lambda(1)\varepsilon(\Lambda), & n=0; \\
               (\lambda\circ S^2)(\Lambda_{(-n)}\Lambda_{(-n-1)}\cdots\Lambda_{(1)}), & n\leq-1.
             \end{array}
           \right.
\end{equation}

The remarkable result due to Y. Kashina, S. Montgomery, S.-H. Ng \cite[Theorem 2.2]{KMN} states that the indicator $\nu_n(A)$ of a finite-dimensional Hopf algebra $A$ is a gauge invariant of $A$ for any $n\geq1$. Later this result has been extended to the case $n\leq0$ by K. Shimizu \cite[Theorem 3.10]{Sh},  where the proof relies on the case of $n\geq1$ and the linear recurrence relation between the Sweedler power maps. In the sequel, we will use Theorem \ref{th2} to give a uniform proof of the two cases.

\begin{thm}
Let $A$ be a finite-dimensional Hopf algebra over the field $\kk$. The $n$-th indicator $\nu_n(A)$ is a gauge invariant of $A$ for any $n\in\mathbb{Z}$.
\end{thm}
\proof
Suppose that $\lambda\in A^*$ is a right integral and $\Lambda\in A$ is a left integral such that $\lambda(\Lambda)=1$. It follows from (\ref{equ01}) that $$\nu_n(A)=(\lambda\circ S)(P_n(\Lambda))\ \ \text{for}\ n\in\mathbb{Z}.$$
Let $A^J$ be the Hopf algebra twisted by a normalized twist $J$ on $A$. The $n$-th indicator of $A^J$ is $\nu_n(A^J)=(\lambda^J\circ S^J)(P^J_n(\Lambda))$. We claim that $\nu_n(A^J)=\nu_n(A)$ for any $n\in\mathbb{Z}$. Indeed,
\begin{align*}
\nu_n(A^J)&=(\lambda^J\circ S^J)(P^J_n(\Lambda))\\
&=\lambda^J(Q_J^{-1}S(P_n^J(\Lambda))Q_J)\\
&=\lambda(\alpha(J^{(1)})S^2(J^{(2)})S(Q_J^{-1})Q_JQ_J^{-1}S(P_n^J(\Lambda))Q_J)\\
&=\lambda(\alpha(J^{(1)})S^2(J^{(2)})S(Q_J^{-1})S(Q_J^{-1}P_n(\Lambda)\alpha(J^{(-1)})S(J^{(-2)})Q_J)Q_J)\ \ \text{by\ Theorem}\ \ref{cor03} (2) \\
&=\lambda(S(P_n(\Lambda))S(Q^{-1}_J)Q_J)\\
&=\lambda(S(P_n(\Lambda))S(Q^{-1}_J)S(J^{(1)})J^{(2)})\\
&=\lambda(\alpha(J^{(2)}_{(1)})S^2(J^{(2)}_{(2)})S(P_n(\Lambda))S(Q^{-1}_J)S(J^{(1)}))\ \ \text{by}\ (\ref{equ4004})\\
&=(\lambda\circ S)(J^{(1)}Q^{-1}_JP_n(\Lambda)\alpha(J^{(2)}_{(1)})S(J^{(2)}_{(2)}))\\
&=(\lambda\circ S)(J^{(1)}Q^{-1}_JP_{n-1}(\Lambda_{(1)})\Lambda_{(2)}\alpha(J^{(2)}_{(1)})S(J^{(2)}_{(2)}))\\
&=(\lambda\circ S)(J^{(1)}Q^{-1}_JP_{n-1}(\Lambda_{(1)}J^{(2)})\Lambda_{(2)})\ \ \ \text{by}\ (\ref{equ005})\\
&=(\lambda\circ S)(J^{(1)}J^{-(1)}S(J^{-(2)})\Lambda_{(1)}P_{n-1}(J^{(2)}\Lambda_{(2)}))\\
&=(\lambda\circ S)(J^{(1)}J^{-(1)}\Lambda_{(1)}P_{n-1}(J^{(2)}J^{-(2)}\Lambda_{(2)}))\ \ \ \text{by}\ (\ref{equ004})\\
&=(\lambda\circ S)(\Lambda_{(1)}P_{n-1}(\Lambda_{(2)}))\\
&=(\lambda\circ S)(P_{n}(\Lambda))\\
&=\nu_n(A).
\end{align*}
If $A'$ is a finite-dimensional Hopf algebra over the field $\kk$ which is gauge equivalent to $A$, then there exists a normalized twist $J$ of $A$ such that $A'$ is isomorphic via a map $\sigma$ to $A^J$ as Hopf algebras. Then $$\sigma\circ S'\circ P'_n=S^J\circ P_{n}^J\circ\sigma,$$ where $S'$ and $P'_n$ are the antipode and the $n$-th Sweedler power map of $A'$ respectively. For any $n\in\mathbb{Z}$, we have
\begin{align*}\nu_n(A')&=\text{tr}(S'\circ P'_{n-1})=\text{tr}(\sigma^{-1}\circ S^J\circ P_{n-1}^J\circ\sigma)\\
&=\text{tr}(S^J\circ P_{n-1}^J)=\nu_n(A^J)=\nu_n(A).
\end{align*}
The proof is completed.
\qed

We need to point out that Theorem \ref{th2} can be used to prove the equality $$\text{tr}(S^n)=\text{tr}((S^J)^n)\ \text{for}\ n=\pm1,\pm2.$$ However, we are not sure Theorem \ref{th2} can be used to prove the equality for all $n\in\mathbb{Z}$. It is known that if $A$ has the Chevalley property, the equality holds for all $n\in\mathbb{Z}$ (see \cite[Theorem 4.3]{CS}).

Finally, we will use Theorem \ref{th2} to show that the Killing
form of a Hopf algebra is invariant under twisting.
Let $A$ be a finite-dimensional Hopf
algebra with antipode $S$ over the field $\kk$. The (left) adjoint representation of $A$ is the map $$\text{ad}: A\rightarrow\text{End}_{\kk}(A),\ a\mapsto\text{ad}a,$$ where the $\kk$-linear map $\text{ad}a:A\longrightarrow A$ is given by
$$(\text{ad}a)(b)=a_{(1)} bS(a_{(2)})\ \text{for}\ b\in A.$$ The Killing form of the Hopf algebra $A$ is defined by
$$(a,b)=\text{tr}(\text{ad}a\circ \text{ad}b)=\text{tr}(\text{ad}ab)\ \text{for}\ a,b\in A.$$
The subspace $$A^\perp =\{ a\in A\mid(a,b)=0\ \text{for\ all}\ b\in A\}$$ is called the Killing
radical of $A.$ This radical is a two-sided ideal of $A$ but not a Hopf ideal in general.

We denote the adjoint representation of the twisted Hopf algebra $A^J$ by $\text{ad}^J$ and the Killing form of $A^J$ by $(-,-)^J$. Namely, $(a,b)^J=\text{tr}(\text{ad}^J(ab)).$
The following result has been given in \cite[Corollary 4.3]{WCL}, we present it here with an alternative proof.

\begin{prop}The Killing form of a finite-dimensional Hopf algebra $A$ is invariant under twisting. Namely,
if $A^J$ is the Hopf algebra twisted by a normalized twist $J$ on a finite-dimensional Hopf algebra $A$, then $(a,b)=(a,b)^J$ for all $a,b\in A$. In particular, $A^{\perp}=(A ^J)^{\perp}$.
\end{prop}
\proof For $a,b\in A,$ we only need to prove that $(a,1)^J=(a,1)$ since $(a,b)^J=(ab,1)^J$ and $(a,b)=(ab,1).$ Suppose that $\lambda\in A^*$ is a right integral and $\Lambda\in A$ is a left integral such that $\lambda(\Lambda)=1$. Recall from Lemma \ref{prop1} that
$$\Delta^J(\Lambda)=\Lambda_{\langle1\rangle}\otimes\Lambda_{\langle2\rangle}=Q_J^{-1}\Lambda_{(1)}\otimes\Lambda_{(2)}S(R)Q_J.$$
By Radford's trace formula (\ref{equ400444}), we have 
\begin{align*}
(a,1)^J&=\text{tr}(\text{ad}^J(a))\\
&=\lambda^J(S^J(\Lambda_{\langle2\rangle})\text{ad}^J(a)(\Lambda_{\langle1\rangle}))\\
&=\lambda^J(S^J(\Lambda_{(2)}S(R)Q_J)\text{ad}^J(a)(Q^{-1}_J\Lambda_{(1)}))\\
&=\lambda^J(Q^{-1}_JS(\Lambda_{(2)}S(R)Q_J)Q_Ja_{\langle1\rangle}Q^{-1}_J\Lambda_{(1)}S^J(a_{\langle2\rangle}))\\
&=\lambda^J(Q^{-1}_JS(\Lambda_{(2)}S(R)Q_J)Q_JJ^{-(1)}a_{(1)}J^{(1)}Q^{-1}_J\Lambda_{(1)}Q^{-1}_JS(J^{-(2)}a_{(2)}J^{(2)})Q_J)\\
&=\lambda(S^2(R^{-1})S(Q^{-1}_J)Q_JQ^{-1}_JS(\Lambda_{(2)}S(R)Q_J)Q_JJ^{-(1)}a_{(1)}J^{(1)}Q^{-1}_J\Lambda_{(1)}Q^{-1}_JS(J^{-(2)}a_{(2)}J^{(2)})Q_J)\\
&=\lambda(S(\Lambda_{(2)})Q_JJ^{-(1)}a_{(1)}J^{(1)}Q^{-1}_J\Lambda_{(1)}Q^{-1}_JS(J^{-(2)}a_{(2)}J^{(2)})Q_J)\\
&=\lambda(S(\Lambda_{(2)})J^{(1)}Q^{-1}_JQ_JJ^{-(1)}a_{(1)}\Lambda_{(1)}Q^{-1}_JS(J^{(2)})S(a_{(2)})S(J^{-(2)})Q_J)\ \ \ \text{by}\ (\ref{equ004})\\
&=\lambda(S^2(S(J^{-(2)})Q_J\leftharpoonup\alpha)S(\Lambda_{(2)})J^{(1)}J^{-(1)}a_{(1)}\Lambda_{(1)}Q^{-1}_JS(J^{(2)})S(a_{(2)}))\ \ \ \text{by}\ (\ref{equ4004})\\
&=\lambda(S(\Lambda_{(2)})J^{(1)}J^{-(1)}a_{(1)}\Lambda_{(1)}S(J^{-(2)})Q_JQ^{-1}_JS(J^{(2)})S(a_{(2)}))\ \ \ \text{by}\ (\ref{equ005})\\
&=\lambda(S(\Lambda_{(2)})a_{(1)}\Lambda_{(1)}S(a_{(2)}))\\
&=\text{tr}(\text{ad}(a))\\
&=(a,1)
\end{align*}
Now $A^{\perp}=(A ^J)^{\perp}$ follows from $(a,b)=(a,b)^J$ for all $a,b\in A$. \qed

As a consequence, we obtain a gauge invariant of $A$ as follows:
\begin{cor}Let $A$ be a finite-dimensional Hopf
algebra over the field $\kk$. Then $\text{dim}_{\kk}(A^{\perp})$ is a gauge invariant of $A$.
\end{cor}

\begin{exa}
Consider the dihedral group $D_n =\{ 1,a,\cdots
,a^{n-1},b,ba,\cdots ,ba^{n-1}\}$ with $a^n =1,\ b^2=1$ and $aba=b$. Suppose $2n\neq0$ in $\kk$. If $n$ is odd, then $(\kk D_n)^{\perp}=0$. If $n$ is even, then  $(\kk D_n)^{\perp}$ is the ideal of $\kk D_n$ generated by $1-a^{\frac{n}{2}}$. This is a Hopf ideal of $\kk D_n$, $\dim_{\kk}((\kk D_n)^{\perp})=n$ and the quotient Hopf algebra is  $\kk D_n/(\kk D_n)^{\perp}\cong\kk D_{\frac{n}{2}}$ (see \cite[Example 3.7]{WCL}).
\end{exa}

\begin{exa}
Let $H_{n,d}$ denote a generalized Taft algebra over the field $\kk$ with $nd\neq0$ in $\kk$. As an algebra, $H_{n,d}$ is generated by $g$ and $h$ subject to the following relations:
$$g^n=1,\  h^d=0,\ hg=qgh,$$
where $d$ divides $n$ and $q$ is a primitive $d$-th root of unity. As a Hopf algebra, the comultiplication $\Delta$, counit $\varepsilon$ and the antipode $S$ of
$H_{n,d}$ are given respectively by
$$\Delta(g)=g\otimes g,\ \Delta(h)=1\otimes h+h\otimes g,\ \varepsilon(g)=1,\ \varepsilon(h)=0,$$
$$S(g)=g^{-1},\ S(h)=-q^{-1}g^{n-1}h.$$ The Hopf algebra $H_{n,d}$ has a $\kk$-basis $\{g^ih^j\mid 0\leq i \leq n-1,0\leq j\leq d-1\}$ and $\dim_{\kk}(H_{n,d})=nd$. If $d=n$, then $H_{n,n}$ is the Taft algebra.
By a straightforward computation, the Killing radical of $H_{n,d}$  is
the ideal of $H_{n,d}$ generated by $g^{d}-1$ and $h$. The dimension of the Killing radical $H_{n,d}^{\perp}$ is $\dim_{\kk}(H_{n,d}^{\perp})=(n-1)d$, which is a gauge invariant of $H_{n,d}$.
\end{exa}

\section*{Acknowledgement}
The first author was supported by Qing Lan Project. The second author was supported by National
Natural Science Foundation of China (Grant No. 11722016). The third author was supported by National
Natural Science Foundation of China (Grant No. 11871063).

\end{document}